%% file: root.tex
\documentclass[letterpaper, 10 pt, conference]{ieeeconf}  % Comment this line out if you need a4paper
\usepackage[T1]{fontenc}

\IEEEoverridecommandlockouts                              % This command is only needed if 
\usepackage{cite}
\usepackage{amsmath,amssymb,amsfonts}
\usepackage{algorithmic}
\usepackage{graphicx}
\usepackage{textcomp}
\usepackage[dvipsnames]{xcolor}
\usepackage{soul}
\usepackage{multirow}

\let\labelindent\relax
\usepackage{enumerate}
\usepackage{enumitem}
\usepackage{comment}
\usepackage{array}
\usepackage{subfigure}
\usepackage{makecell,tabularx, float}

\makeatletter
\let\NAT@parse\undefined
\makeatother
\usepackage[pagebackref=true]{hyperref}
 \hypersetup{ %
  colorlinks=true, %
  pdfstartview={FitH}, %
  citecolor=Black, %
  linkcolor=Black, %
  urlcolor=Blue %
 }

\title{\LARGE \bf
A Multi-objective Mixed-integer Programming Approach for Supply Chain Disruption Response with Lead-Time Awareness*
}

\author{Juan-Alberto Estrada-Garcia$^{1}$, Mingjie Bi$^{2}$, Dawn M. Tilbury$^{2,3}$, Kira Barton$^{2,3}$ and Siqian Shen$^{1}$% <-this % stops a space
\thanks{*This work was funded in part by The United States National Science Foundation (NSF) grant \#CMMI-2034974.}% <-this % stops a space
\thanks{$^{1}$Juan-Alberto Estrada-Garcia and Siqian Shen are with the Department of Industrial and Operations Engineering, 
        University of Michigan, Ann Arbor, MI 48109, USA,
        {\tt\small \{juanest, siqian\}@umich.edu}}%
\thanks{$^{2}$Mingjie Bi is with the Robotics Department, 
        University of Michigan, Ann Arbor, MI 48109, USA,
        {\tt\small mingjieb@umich.edu}}%
\thanks{$^{2,3}$Dawn M. Tilbury and Kira Barton are with the Robotics Department and the Department of Mechanical Engineering, 
        University of Michigan, Ann Arbor, MI 48109, USA,
        {\tt\small \{tilbury, bartonkl\}@umich.edu}}%
}

\begin{document}

\allowdisplaybreaks
\maketitle
%\thispagestyle{empty}
%\pagestyle{empty}

%%%%%%%%%%%%%%%%%%%%%%%%%%%%%%%%%%%%%%%%%%%%%%%%%%%%%%%%%%%%%%%%%%%%%%%%%%%%%%%%
\begin{abstract}

Supply chain (SC) risk management is influenced by both spatial and temporal attributes of different entities (suppliers, retailers, and customers). Each entity has given capacity and lead time for processing and transporting products to downstream entities. Under disruptive events, lead time and capacities may vary, which affects the overall SC performance. There have been many studies on SC disruption mitigation, but often without considering lead time and the magnitude of lateness. In this paper, we formulate a mixed-integer programming (MIP) model to optimize SC operations via a routing and scheduling approach, to model the delivery time of products at different entities as they flow throughout the SC network. We minimize a weighted sum of multiple objectives involving costs related to transportation, shortage, and delivery lateness. We also develop a discrete-event simulation framework to evaluate the performance of solutions to the MIP model under lead time uncertainty. Via extensive numerical studies, we show how the attributes of SC entities affect the performance, so that we can improve SC design and operations under various uncertainties.


\end{abstract}

\section{Introduction}
\label{sec:introduction}

%%%%%%%%%% Motivation & Problem statement %%%%%%%%%%
Supply chain (SC) design and operations are impacted by the capacities and production/processing lead time at different entities in an SC network, as well as how the entities are connected and communicate. %~\cite{mentzer2001defining}. 
The global modern SC network has become more complex, leading to more frequent disruptions that affect SC performance. %~\cite{xu2020disruption}. 
In addition to capacities, the lead time of procurement and production faces increasing uncertainties due to labor shortages and unpredictable global environments in recent years~\cite{hammami2013optimisation}.

To cope with lead-time disruptions, enterprises need to incorporate both spatial and temporal attributes of flows and productions into SC design and management. %~\cite{park2010three}.
It is crucial to understand how lead time impacts disruption response and SC risk management~\cite{chang2019effect}. 
In this work, we consider the following decision-making problem: Given a layout of an SC network, existing flows between different entities, and disrupted capacities and/or lead time, how do shortage and lateness penalties, as well as different attributes of the SC entities, affect response decisions?

%%%%%%%%%% Problem statement & Lit review %%%%%%%%%%
The existing literature on the optimization of SC networks' operations models lead time either as a constraint or as an element of the objective function.
In~\cite{hammami2013optimisation}, the authors use a mixed-integer programming (MIP) model for multi-echelon SC considering lead time as hard constraints.
They perform extensive numerical experiments to study the feasibility of a network-flow problem with hard lead-time requirements.
In~\cite{you2008design}, authors propose a multi-period, mixed-integer nonlinear programming (MINLP) model to optimize SC design. They minimize the expected lead time of products, by considering the expected lead time of flows. However, these approaches do not consider lead time disruptions and the corresponding SC response. They also do not model penalties for delivery lateness. The existing literature often models lead time only for end customers, instead of each entity in the SC~\cite{petridis2015optimal}. For example,~\cite{eskigun2005outbound} introduces a capacitated network design model but only considers the transportation lead time and penalizes the late delivery of products to the final clients given a homogeneous fixed penalty; 
\cite{aqlan2016supply} poses a multi-objective goal-programming approach to optimize the planning of multi-period SCs given non-negligible lead time; \cite{hahn2012value} introduces a two-stage multi-period stochastic optimization model with a scenario-based solution approach, to enhance the robustness of SC design. However, the formulations in both~\cite{aqlan2016supply} and~\cite{hahn2012value} only model delivery lateness as back-orders that always ship in the next period.

All the above methods focus on the initial SC network design and planning to achieve network-level lead-time reduction.
However, for disruption response, production and flows that require to be recovered may exist at any entity in an SC network. 
Enterprises need to determine whether and how to vary the production and flow plans to meet customer demands. 
Therefore, delivery times of products are critical elements to be tracked through the entire SC network, so that penalties for product lateness at any point in the SC network can be considered. To the best of our knowledge, no existing work fully addresses such needs and challenges.

Furthermore, %to make quick disruption-response decisions, 
understanding how different SC attributes impact the performance from both network and entity levels is critical. Most existing work investigates the effects of network typologies on SC resilience, without the consideration of lead time \cite{kim2015supply, nair2011supply, zhao2022supply}. The study in \cite{chang2019effect} examines the lead-time effect on resilience across different stratifications in an SC. However, the authors did not discuss how the SC attributes play a role in the disruption response. 

%%%%%%%%%% Contribution %%%%%%%%%%

To address all the above limitations, the main contributions of this work are as follows. (i) We derive an MIP formulation that tracks  temporal attributes of flows through an SC network, while incorporating both fixed and additive penalties on delivery lateness. (ii) We develop a discrete-event simulation framework to measure the performance of the MIP solutions and the effects of lead-time uncertainty. (iii) We investigate the effects of lateness penalty and network topology on the performance of lead-time disruption response through several comprehensive SC instances.

%%%%%%%%%% Paper structure %%%%%%%%%%
The rest of the paper is organized as follows. In Section \ref{sec:centralizedmodel}, we introduce detailed notation and the multi-objective MIP model for deriving new flows and response solutions to given SC disruptions. In Section \ref{sec:outofsampletest}, we present the discrete-event simulation framework for evaluating the performance of the MIP model under uncertain lead time. In Section \ref{sec:casestudy}, we demonstrate the results of a simulated case study and provide managerial insights. Finally, in Section \ref{sec:conclusion}, we conclude the paper and propose future research directions.

%% Section 2

\section{Multi-objective MIP}
\label{sec:centralizedmodel} 
%We first introduce the mathematical model proposed in~\cite{bi2022model} and extend it to an MIP model that tracks product flows and arrival time at each entity in an SC network, such that we can incorporate the magnitude of all late deliveries in the multi-objective function.

\subsection{Notation and assumptions}
\label{sec:NotationAssumpt}

Consider a directed graph $G(V,E)$ representing an SC network with a vertex set $V$ of all SC entities and an edge set $E$ of all their connections. We denote $V^{\text{d}} \subset V$ as the subset of entities that distribute products, named distributors (i.e., no transformation of components occurs within them); $V^{\text{o}} \subset V$ denotes the subset of entities where transformations of products occur, named OEMs; $V^{\text{s}} \subset V$ is the subset of entities who supply products and raw materials, named suppliers (i.e., having no upstream flows). Edges in set $E$ convey potential flows of products and resources between entities. We denote $K$ as the set of all products and component types within the SC. Additionally, for each $k \in K$, we have a subset $K'(k) \subset K$ of components $k'$ required for the production of $k$.
We make the following assumptions in this paper: %{\color{red} (i) At the beginning of the period, all entities receive demand information to fulfill. (ii) If entity $j$ orders product $k$ from entity $i$, all the flows will be sent simultaneously. (iii) Suppliers begin production immediately. (iv) OEMs will start producing/sending product $k$ until all required upstream products have been received. (5) Distributors wait until all units of product $k$ are collected.} 
(i) Knowledge of disruptions and response time are immediate.  (ii) Flow of product $k$ from $i$ to $j$ is treated as an indivisible unit. (iii) All entities wait until all required upstream flows have been received before sending their downstream flows. We summarize the parameters and decision variables used in this work in~Table~\ref{tab:notation}.  

\subsection{Multi-objective MIP for SC disruption response with lead time awareness}
\label{subsec:MIPmodel}

\input{tables/notations}

We compute the objective of our MIP as the total SC operation cost in \eqref{eq:obj}, where $y, x, I, p, \beta, \zeta, \Delta, a, o, z,w $ are the decision variables shown in Table \ref{tab:notation}. We model the cost of transportation, manufacturing, and product holding, respectively in \eqref{obj:var}, fixed costs for transportation route setting, manufacturing capacity in~\eqref{obj:fix}, and the penalties for unmet demand and late deliveries in \eqref{obj:penalties}.
\begin{subequations}\label{eq:obj}
\begin{align}\small
\mathcal{J}&(y,\beta,x,p,\zeta,I,\Delta,a,o,z,w) \nonumber\\
=&\sum_{(i,j) \in E, k \in K} c_{ijk} y_{ijk} + \sum_{i \in V, k \in K} h_{ik}I_{ik} + \sum_{i \in V, k \in K} e_{ik}p_{ik} \label{obj:var}\\
&+\sum_{(i,j) \in E, \ k \in K} f_{ij}\beta_{ijk}+
\sum_{i \in V} \phi_{i} \zeta_{i} \label{obj:fix}\\
%\sum_{i \in V} \eta_{i} \xi_{i} + \label{obj:fix} \\
&+\sum_{i \in V, k \in K}\rho_{ik}^d \Delta_{ik}^d + \sum_{(i,j) \in E, k \in K} w_{ijk}, \label{obj:penalties}
% \sum_{i \in V,\ k \in K} \rho_{ik}^s\Delta_{ik}^I +\label{obj:pe}\\
%\sum_{(i,j) \in E}\rho_{ij}^{E}\Delta_{ij}^{E} + \sum_{i \in V}\rho_i^{V}\Delta_i^V \label{obj:re}
\end{align}
\end{subequations}
\normalsize
We present the overall MIP model to optimize SC network flows and operations in~\eqref{mip}:
\begin{subequations}\label{mip}
\begin{align}\small
\label{obj:mip}
\min_{X}& \hspace{1ex}\mathcal{J}\\
%\sum_{(i,j) \in E}\left(f_{ij}z_{ij}+\rho_{ij}^{E}\Delta_{ij}^{E} \right) + \sum_{i \in V}\left(\phi_{i} \zeta_{i} + \eta_{i} \xi_{i} + \rho_i^{V}\Delta_i^V\right) +\nonumber\\
%&\sum_{k \in K}\left(
%\sum_{(i,j) \in E}c_{ijk}y_{ijk} + \sum_{i \in V}\left(h_{ik} I_{ik} + e_{ik}p_{ik} + \rho_{ik}^d \Delta_{ik}^d \right)\right) \label{eq:obj}\\
%  + \rho_{ik}^s\Delta_{ik}^I
\text{s.t.}
&\label{cst:flowbalance}
\sum_{j: (i, j)\in E}y_{ijk} - \sum_{j: (j, i)\in E}y_{jik} +  \sum_{k' \in K}r_{k'k}p_{ik'}  \nonumber\\
& \hspace{0ex}- p_{ik} = x_{ik} + I_{ik}^0 - I_{ik}, \ \forall \ i \in V,\ k \in K, \\
& y_{ijk} \leq q_{ij} \beta_{ijk}, \ \forall \ (i,j) \in E, k \in K, \label{cst:individualcap} \\
& \sum_{k \in K} y_{ijk} \leq q_{ij}, \ \forall \ (i,j) \in E \label{cst:mixcap}, \\
  \label{cst:cap_p}
& \sum_{k \in K}p_{ik} \leq \bar{p}_{i}\zeta_{i},\ \forall \ i \in V,\\
%   \label{cst:cap_I}
%& \sum_{k \in K}I_{ik} \leq \bar{I}_{i}\xi_{i},\ \forall \ i \in V\\
\label{cst:pnlty_d}
& \Delta_{ik}^d \geq x_{ik}- d_{ik},\ \forall \ i \in V,\ k \in K,\\
%   \label{cst:pnlty_I}
% & \Delta_{ik}^I \geq I_{ik}^s - I_{ik} ,\ \forall i \in V,\ k \in K\\
%& \Delta_i^E \geq z_{ij} - z_{ij}^{E_0}, \forall (i, j) \in E\\
%& \Delta_i^E \geq z_{ij}^{E_0} - z_{ij}, \forall (i, j) \in E\\
%& \Delta_i^V \geq \zeta_i - \zeta_i^{V_0}, \forall i \in V\\
%& \Delta_i^V \geq \zeta_i^{V_0} - \zeta_i, \forall i \in V\\
%&  y_{ijk}, I_{ik}, L_{ik}, \Delta_{ik}^d \geq 0,\zeta_{i}, z_{ij}, \xi_{i}, \Delta_i^V, \Delta_{ij}^E  \in \{0, 1\}, \nonumber\\
&a_{ijk} = (l_{ijk}+o_{ik})\beta_{ijk},\forall \ (i,j) \in E,k \in K , \label{cst:acomp}\\
    &o_{jk} \geq a_{ijk'}, \nonumber \\& \ \forall \ (i,j) \in E,i \in V \setminus V^{\text{s}},\ k \in K,k'\in K'(j,k), \label{cst:otransf}\\
    &o_{ik} = 0, \ \forall \ i \in V^{\text{s}},k \in K, \label{cst:opart}\\
    &a_{ijk} \leq t_{jk} + \mathcal{M}z_{ijk},\forall \ (i,j) \in E,k \in K, \label{cst:zbigM}\\
    &w_{ijk} = z_{ijk}(c^{\text{f}}_{ijk}+(c^{\text{u}}_{ijk}(a_{ijk}-t_{jk}))), \nonumber \\
    & \forall \ (i,j) \in E,k \in K , \label{cst:penaltyw}\\
&y_{ijk},x_{ik},I_{ik},\Delta_{ik}^{d},a_{ijk},o_{ik},w_{ijk} \geq 0,\nonumber \\
    &\forall i \in V,(i,j) \in E, k \in K,\label{cst:yxIdeltadom}\\
     &\zeta_{i},\beta_{ijk},z_{ijk} \ \in \{0,1\}, \nonumber \\
     & \forall \ i \in V, \ \forall \ (i,j) \in E, k \in K, \label{cst:zetabetadom}
\end{align}
%The objective function \eqref{eq:obj} rearranges the original objective \eqref{obj-mip}. 
\end{subequations}\normalsize
where we denote the vector of decisions as $X = [y, x, I, p, \beta, \zeta, \Delta, a, o, z, w]^\top $. 
The constraints include: 
\subsubsection{Flow balance}
Constraints~\eqref{cst:flowbalance} balance the flow of products at each entity in the SC. 
\subsubsection{Capacities}
In~\eqref{cst:individualcap},~\eqref{cst:mixcap},~\eqref{cst:cap_p}, flows on each edge and productions at each entity are restricted by given capacities.
\subsubsection{Unmet demand penalty}
Constraints~\eqref{cst:pnlty_d} compute the unsatisfied demand for products at each entity. 
\subsubsection{Delivery times}
With constraints~\eqref{cst:acomp}, we model the delivery time of products at different entities. 
Note that we can model products with different processing and transportation lead times. With constraints \eqref{cst:otransf}--\eqref{cst:opart} we compute the time at which downstream flows are ready to be processed (i.e., variable $o_{ik}$), depending on the readiness of upstream products the entity requires. 
In constraints~\eqref{cst:zbigM}, we compare the delivery time of products $a_{ijk}$ with the due date $t_{jk}$, such that we can penalize the flow whenever the delivery is late (via a big-M approach). We can introduce heterogeneous delivery deadlines at any entity of the SC, such that we can model SC networks that not only require the on-time delivery of products to final customers but also at intermediate stages of the SC network. Finally, we model the late delivery penalty as a function of lateness in constraints~\eqref{cst:penaltyw}, where we have both fixed penalty cost (for any lateness) and unit penalty (for each unit of lateness for each flow) at each entity in the SC. 
%For example, one can penalize any lateness significantly (i.e., having a large fixed penalty and comparatively small unit penalty), or accept some degree of lateness (i.e., having a small or no fixed penalty). We note that if the penalty coefficients are all set to zero, it reduces to a lead-time-neutral solution to the problem.
\subsubsection{Variable domains}
Constraints in~\eqref{cst:yxIdeltadom} and~\eqref{cst:zetabetadom} specify the domain of decision variables.

%% Out of sample testing
\section{Out-of-sample testing}
\label{sec:outofsampletest}

%In this section, we present the baseline SC network instance we use to perform numerical experiments. Then, we introduce the simulation framework to test the performance of solutions to the MIP model~\eqref{mip} under stochastic lead time. Lastly, we describe the case study scenarios we test the response decisions for different SC network disruptions.

\subsection{SC instance}
\label{sec:SCinstance}

\begin{figure}[tb]
\centerline{\includegraphics[width=0.95\columnwidth]{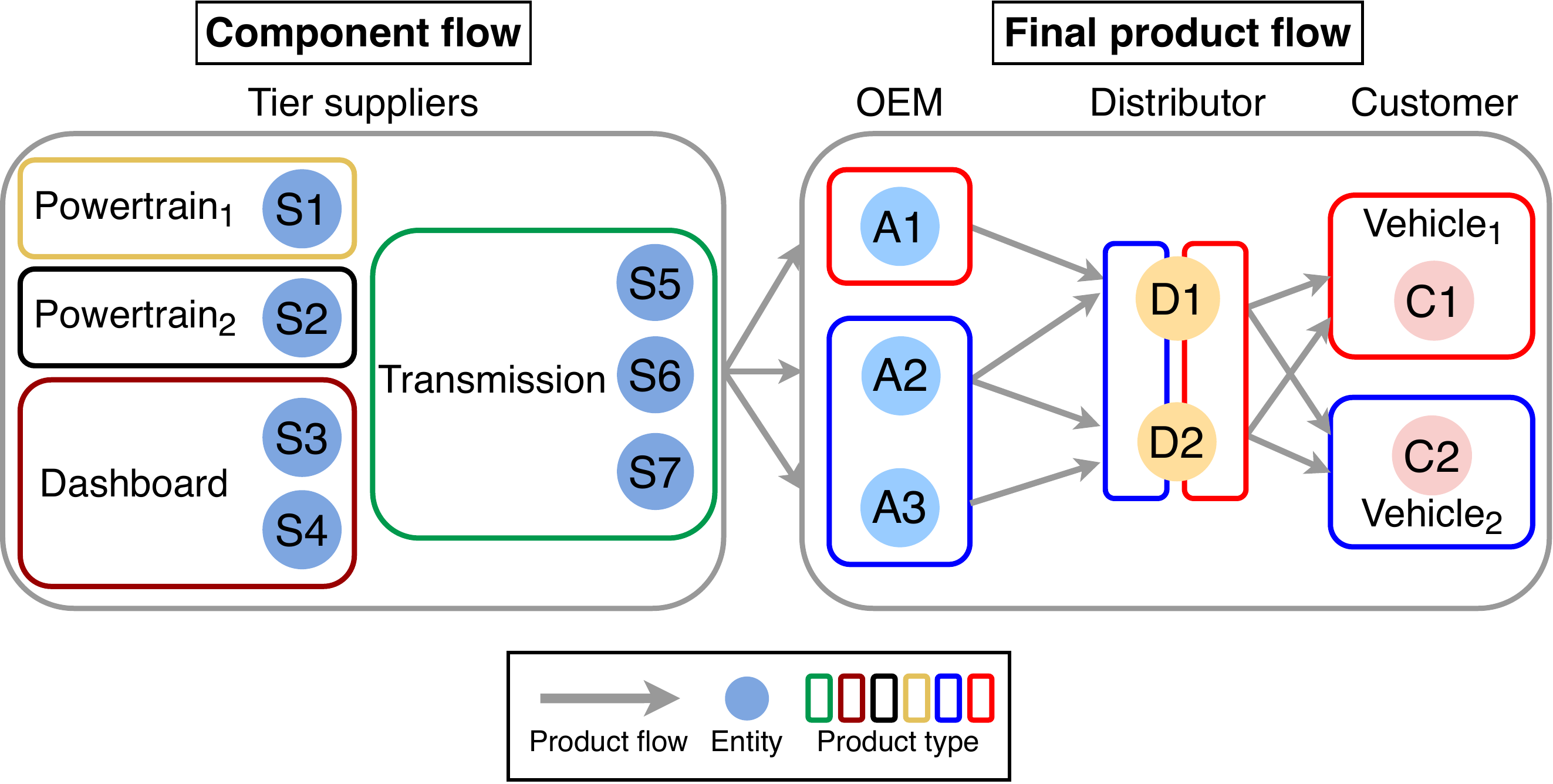}}
\caption{The SC network for the case study. }
\label{fig:setup}
\end{figure}

The baseline SC network is an automotive SC adapted from~\cite{willems2008data}.
%In this work, we illustrate the implementation of the proposed framework with a simulated supply chain network adapted from~\cite{willems2008data}. 
The network consists of 7 suppliers, 3 OEMs, 2 distributors, and 2 end customers. 
The business relationship established between entities indicates the networks' topology, as shown in Fig.\ \ref{fig:setup}.
The annotations represent product types that the suppliers and OEMs can produce and hold in inventory.

\subsection{Simulation framework}
\label{sec:SimulationFramework}

We use discrete-event simulation to evaluate optimal solutions to Model~\eqref{mip} under lead-time uncertainty. 
%We solve the MIP model~\eqref{mip} to compute optimal product flows $y_{ijk}$ considering a deterministic lead time, which can be, e.g., the expected lead time based on historical data. Then, we input $y_{ijk}$ to the SC network $G$ from upstream to downstream, obtaining the resulting arrival time of flows at each entity until all products are delivered to end customers.
%
% We consider the optimal value for variables $y_{ijk}$ as the input flows for the simulation model. Other inputs include the graph $G$ describing the network, the subsets of entity types (i.e., distributor, transformer, or initial supplier), the product structure (i.e., the subsets of subproducts-product relationships), and finally, the lead times $l_{ijk}$ used for the deterministic model optimization.
%
\begin{figure}[tb]
\centerline{\includegraphics[width=\columnwidth]{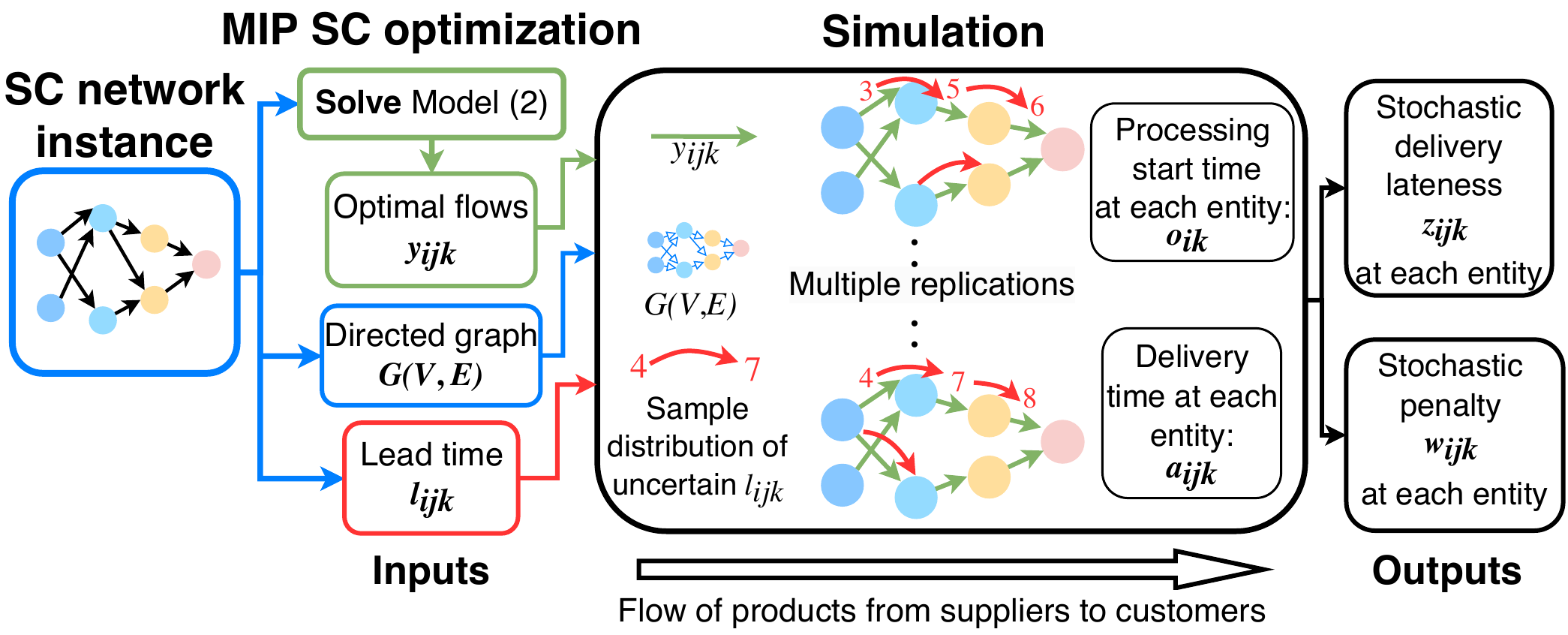}}
\caption{Simulation framework for out-of-sample testing.}
\label{fig:simulation}
\end{figure}
Fig. \ref{fig:simulation} depicts inputs for the simulation model, general simulation process, and simulation outputs. We initialize the simulation with the flows from the suppliers (i.e., entities in subset $V^{s}$). We assume independence in the distribution of lead time for each entity.
The stochastic lead time is sampled from a log-normal probability distribution with $\mu=l_{ijk}$ (i.e., the deterministic lead time $l_{ijk}$ used in Model \eqref{mip}) and $\sigma=0.3$, such that whenever there exits a flow, the realized lead time is a sample from that distribution. We store delivery information at every entity using the value of variable $a_{ijk}$. Having computed this information, we consider the downstream entities that just received flows and use the product structure information to compute their starting time for processing the downstream flows, using the value of variable $o_{ik}$. We continue this process iteratively until the final products are delivered to all end clients. We run the simulation for multiple replications in parallel to analyze the SC performance under different realizations of the lead time.

\subsection{Case study scenario design}
\label{sec:scenarioDesign}

\subsubsection{Baseline scenario}
\label{sec:baselineScenario}
We solve the MIP model~\eqref{mip} for the instance described in Section \ref{sec:SCinstance}, by selecting parameters so that all can deliver products on time. We consider this solution as the baseline benchmark of our case study. 
% \revb{consider the optimal solution to the case study described in Section \ref{sec:SCinstance} as the baseline of our case study, which we compare to other SC network layouts and parameters.} 

 \subsubsection{Delivery lateness policies}
 \label{sec:ArrLatePoli}

We perform a series of studies to analyze the effect of lateness penalties on SC disruption responses under uncertain lead time. We interpret different unit-penalty values as policies chosen by decision makers and model these policies through different ratios between the unit and fixed cost of penalizing lateness.
Instances, in which lateness does not matter, are considered to have neither unit nor fixed penalty. (This configuration is equivalent to the model considered in our prior work \cite{bi2022model} in which we do not study lead time or lateness.) For other instances, we consider a policy that only contains unit penalty (ratio 1:0), and  policies that consider increasing fixed penalties (ratios 1:500 and 1:5000), to model situations in which being late by any amount is undesirable. The values 500 and 5000 as fixed penalty cost are chosen to be proportional to the volume of flows in the current instance such that the fixed penalty dominates the total delivery lateness penalization as the ratio increases.
In practice, these different lateness policies could happen in different industries, e.g., in industries that frequently allow back-orders, industries with perishable products, or those that fit in a larger SC with tight due dates. 

\subsubsection{SC entity disruptions}
\label{sec:SCentityDisrupts}

 We identify the baseline SC network instance (Fig. \ref{fig:baseline_flows}) as a hierarchical network structure.  Sequential connections between different entities exist \cite{kim2015supply}. This configuration is justified by the process of assembling different components to yield products for end customers. 

 We study the effect of entity depth (i.e., how far back in the product flow process they are positioned) on SC disruption response.  We design three scenarios considering the SC network in Fig.\ \ref{fig:setup}: (i) Disrupt a transmission supplier with a depth three by doubling their lead time. (ii) Disrupt an OEM  by doubling their lead time. (iii) Disrupt a distributor of vehicle$_1$ and vehicle$_2$ by tripling their lead time. %We select parameters to model severe disruptions on entities that are relied upon in the optimal solution of the baseline.

\subsubsection{SC topology}
\label{sec:SCtopologyIntro}

The baseline SC network has multiple entities capable to supply, assemble, or distribute materials and final products. Its topology follows a ``tree-like'' structure in which the in-degree of entities reduces as the sequence of entities gets closer to end customers. This strategy has high cost due to multiple contractual relationships with different suppliers, but can provide the flexibility to respond to disruptions by choosing backup entities \cite{nair2011supply}.

We study the resiliency of different network configurations to understand their performance under entity disruption and lead-time uncertainty. We compare the baseline configuration with two others as follows. (i) A ``reverse tree-like'' structure, in which the budget for building redundant relationships is directed to entities closer to end customers (i.e., we remove one transmission supplier and add a fourth OEM and a third distributor). (ii) A ``chain-like'' structure with fewer backup entities, without removing productive capacities in the SC (i.e., we remove one transmission supplier, one OEM, and one distributor). This configuration has the least cost, as the economy of scale benefits from having larger contracts with fewer entities. For all typologies, we have the same total capacity between all the entities and the same demand. We optimize the non-disrupted instance, and then disrupt similar agents in each topology setting, to compare out-of-sample performance of the response given by the original solution.

\section{Case study}
\label{sec:casestudy}

\begin{figure}[tb]
\centerline{\includegraphics[width=0.5\columnwidth]{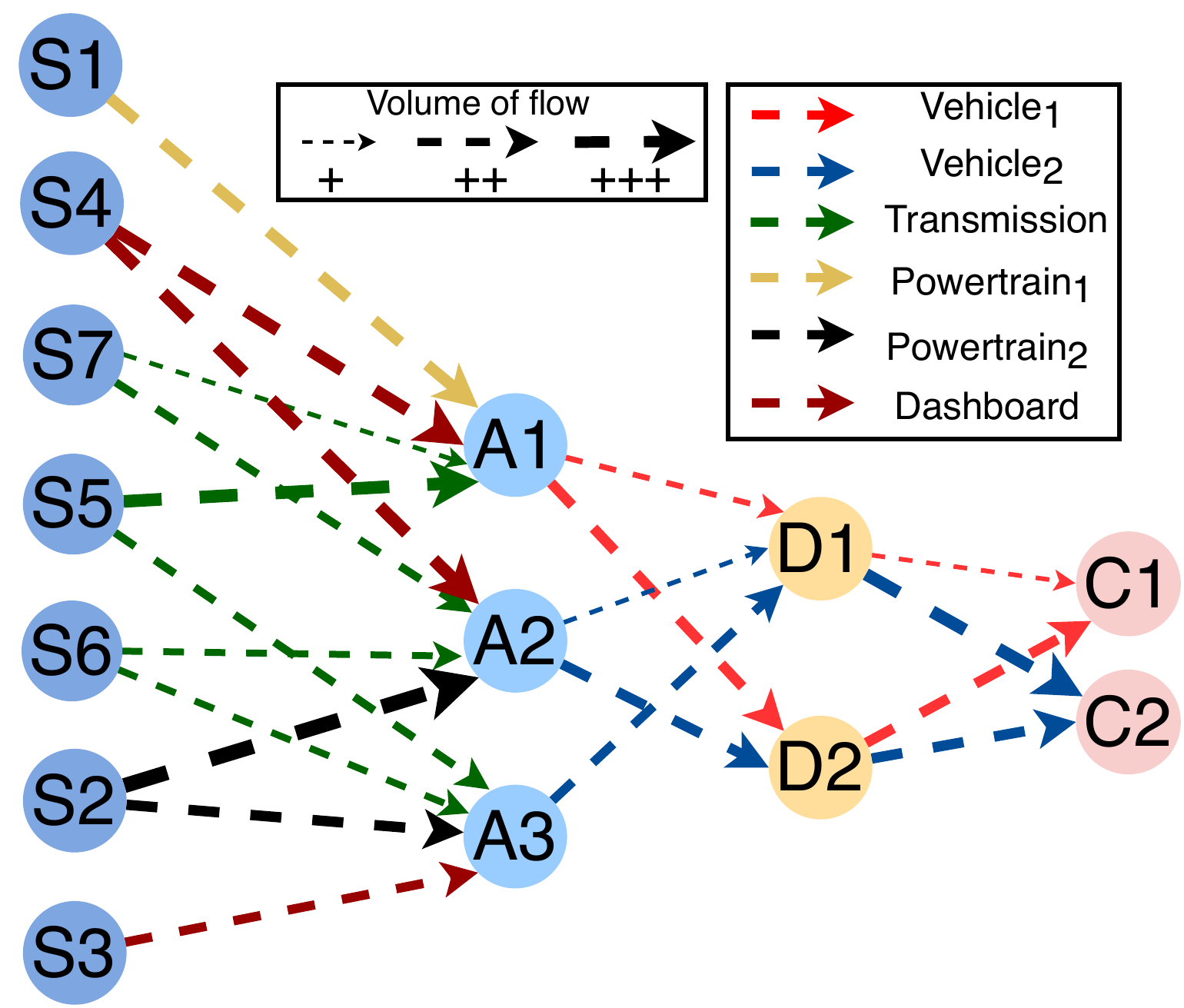}}
\caption{Optimal flows from solving MIP model~\eqref{mip} for baseline instance.}
\label{fig:baseline_flows}
\end{figure}

%
%\begin{figure}[tb]
%\centerline{\includegraphics[width=0.8\columnwidth]{figures/baseline_solution_simulation.pdf}}
%\caption{Out-of-sample performance of MIP model~\eqref{mip} solution to the baseline instance.}
%\label{fig:baseline_simulation}
%\end{figure}

\begin{figure*}[t]
\centering
\subfigure[Disruption at S7.]{\includegraphics[width=0.6\columnwidth]{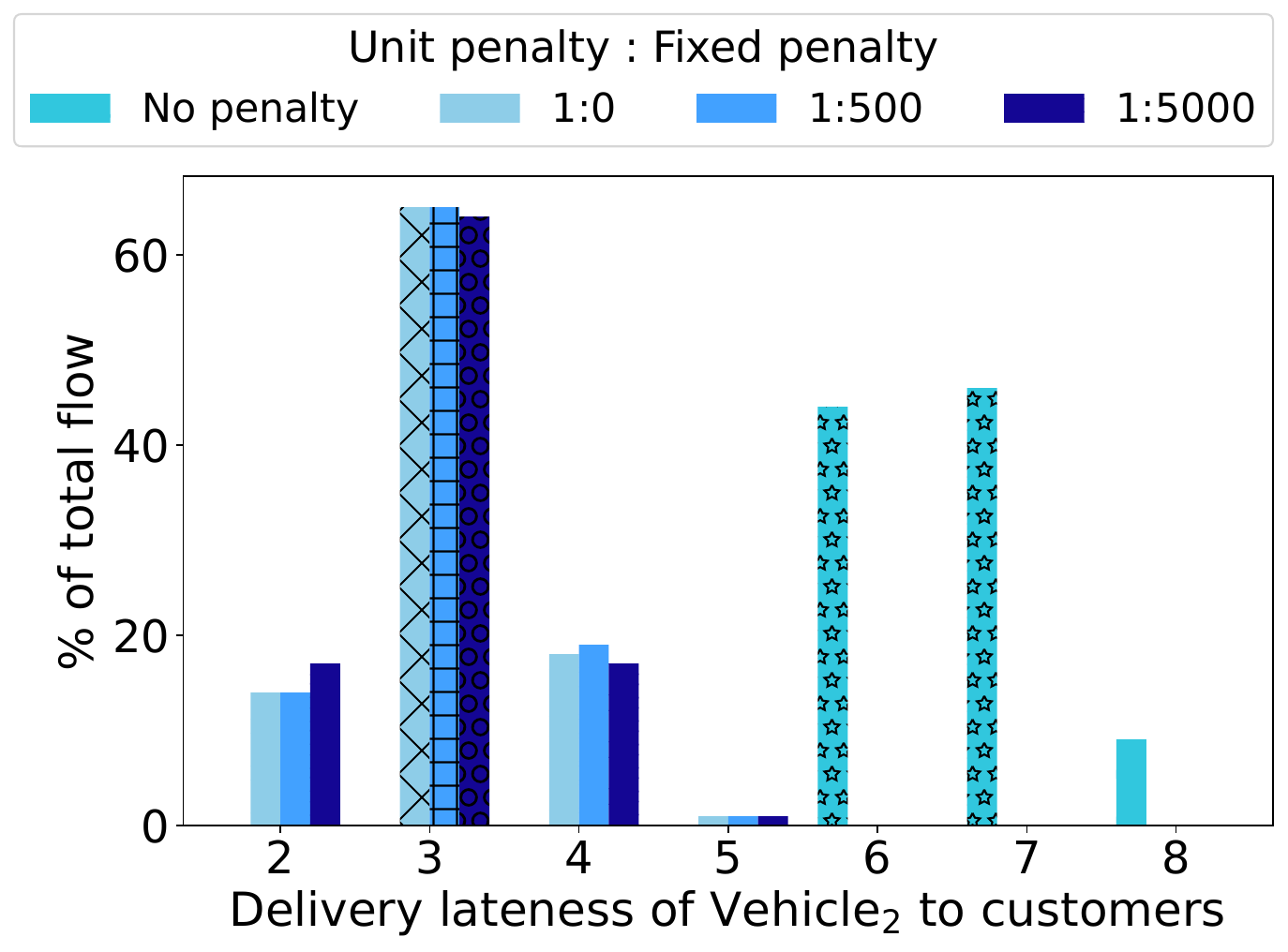}%
\label{fig:S7}}
\hfil
\subfigure[Disruption at  A2.]{\includegraphics[width=0.6\columnwidth]{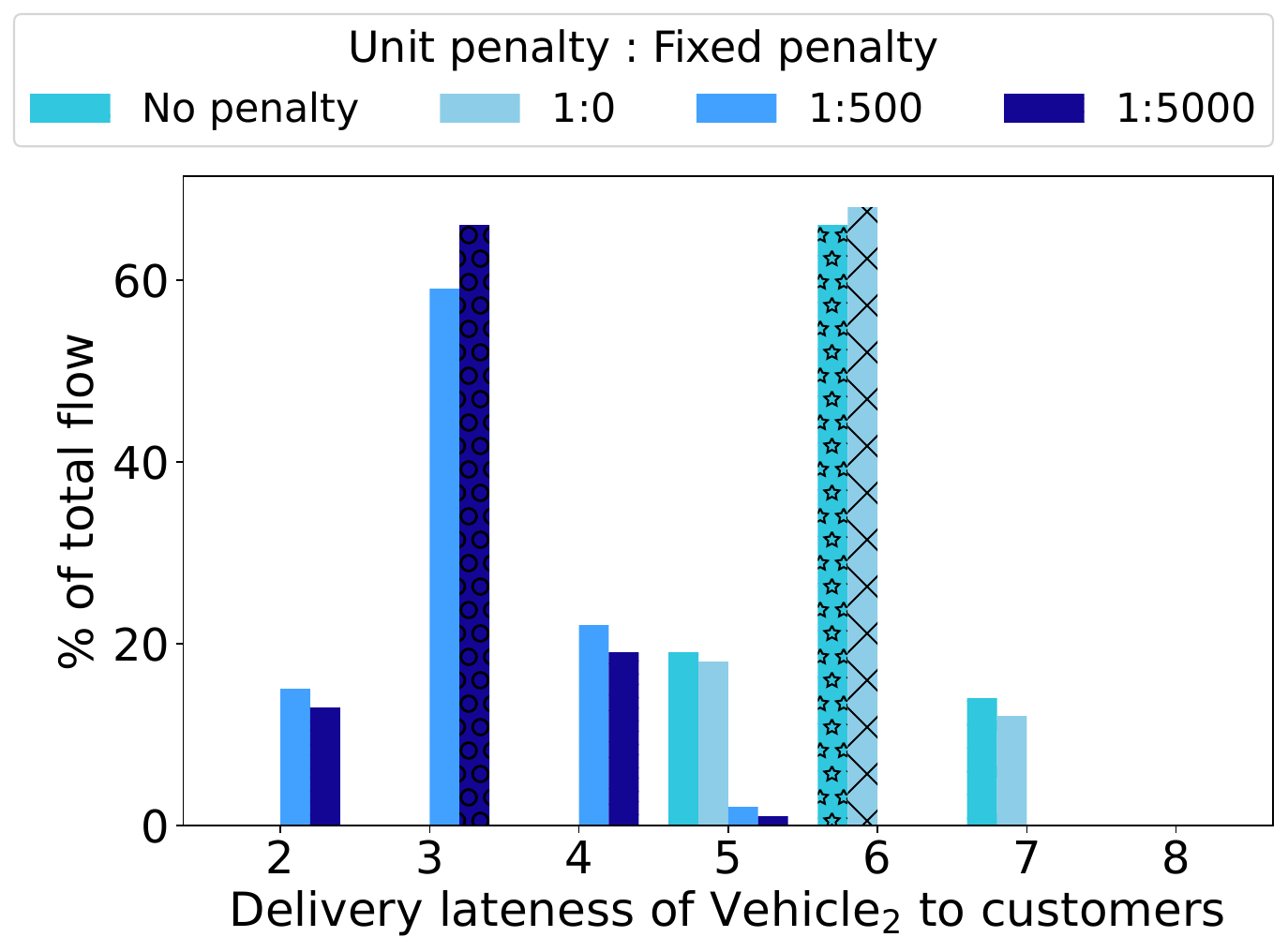}%
\label{fig:A2}}
\hfil
\subfigure[Disruption at D1.]{\includegraphics[width=0.6\columnwidth]{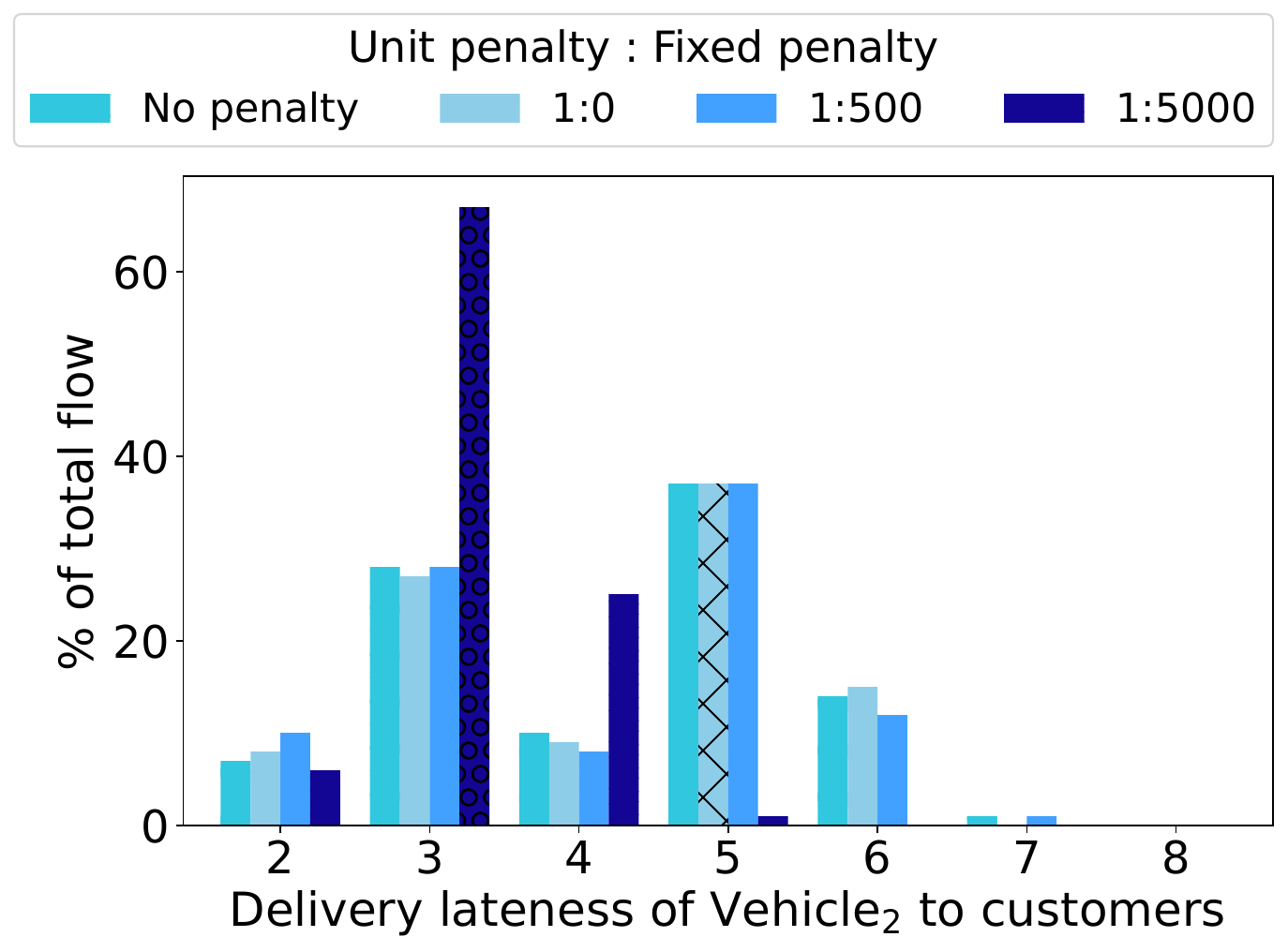}%
\label{fig:D1}}
\caption{Out-of-sample performance of solutions under different lateness policies and entity disruptions for a final product of the SC.}
\label{fig:resultsb_simulation}
\end{figure*}

\begin{figure*}[t]
\centering
\subfigure[Optimal strategy K for no lateness penalty and ratio 1:0]{\includegraphics[width=0.6\columnwidth]{figures/baseline_solutionflows.pdf}%
\label{fig:flows_A2_rat_01}}
\hfil
\subfigure[Optimal strategy E for penalty ratio 1:500]{\includegraphics[width=0.6\columnwidth]{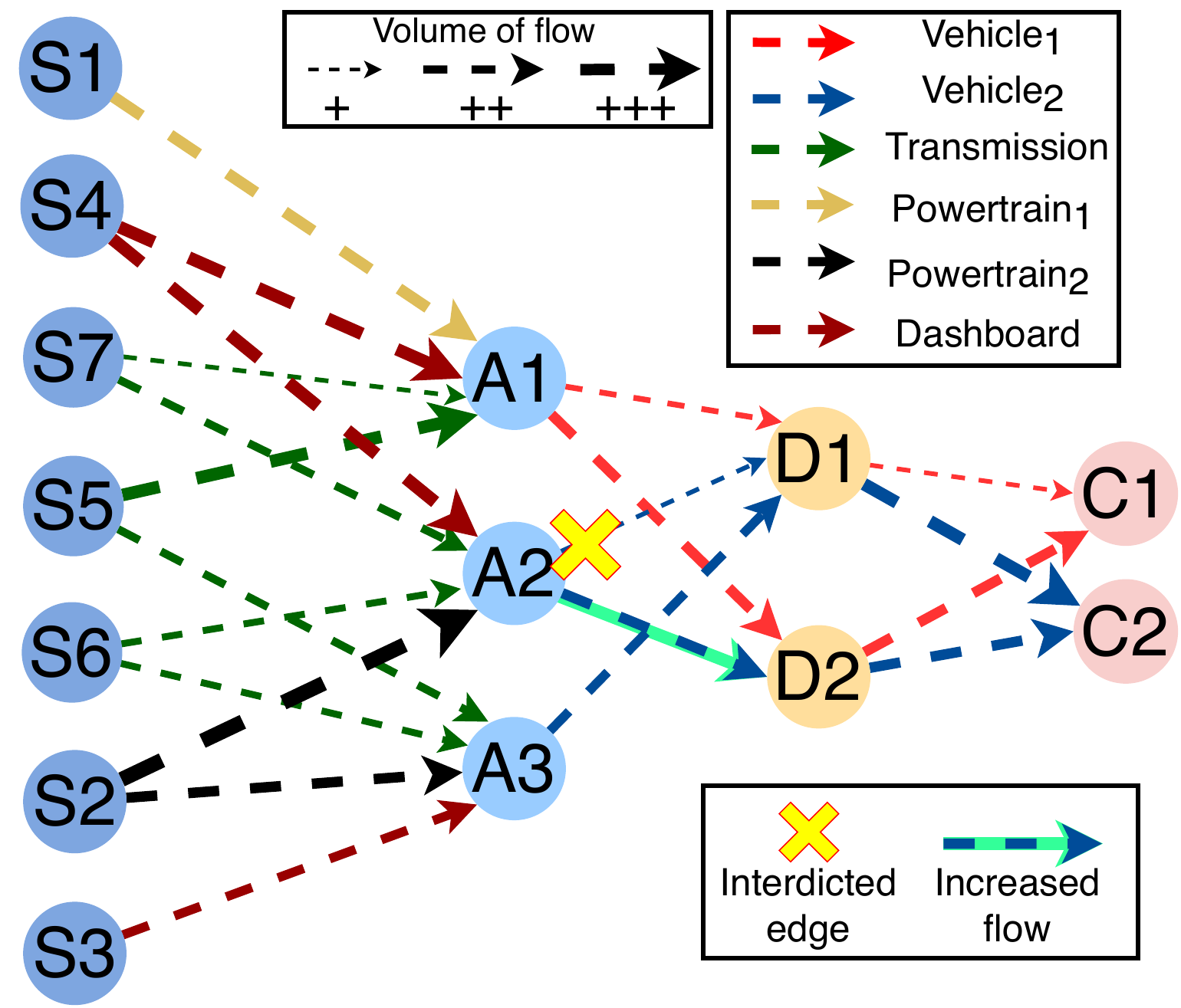}%
\label{fig:flows_A2_rat_500}}
\hfil
\subfigure[Optimal strategy V for penalty ratio 1:5,000]{\includegraphics[width=0.6\columnwidth]{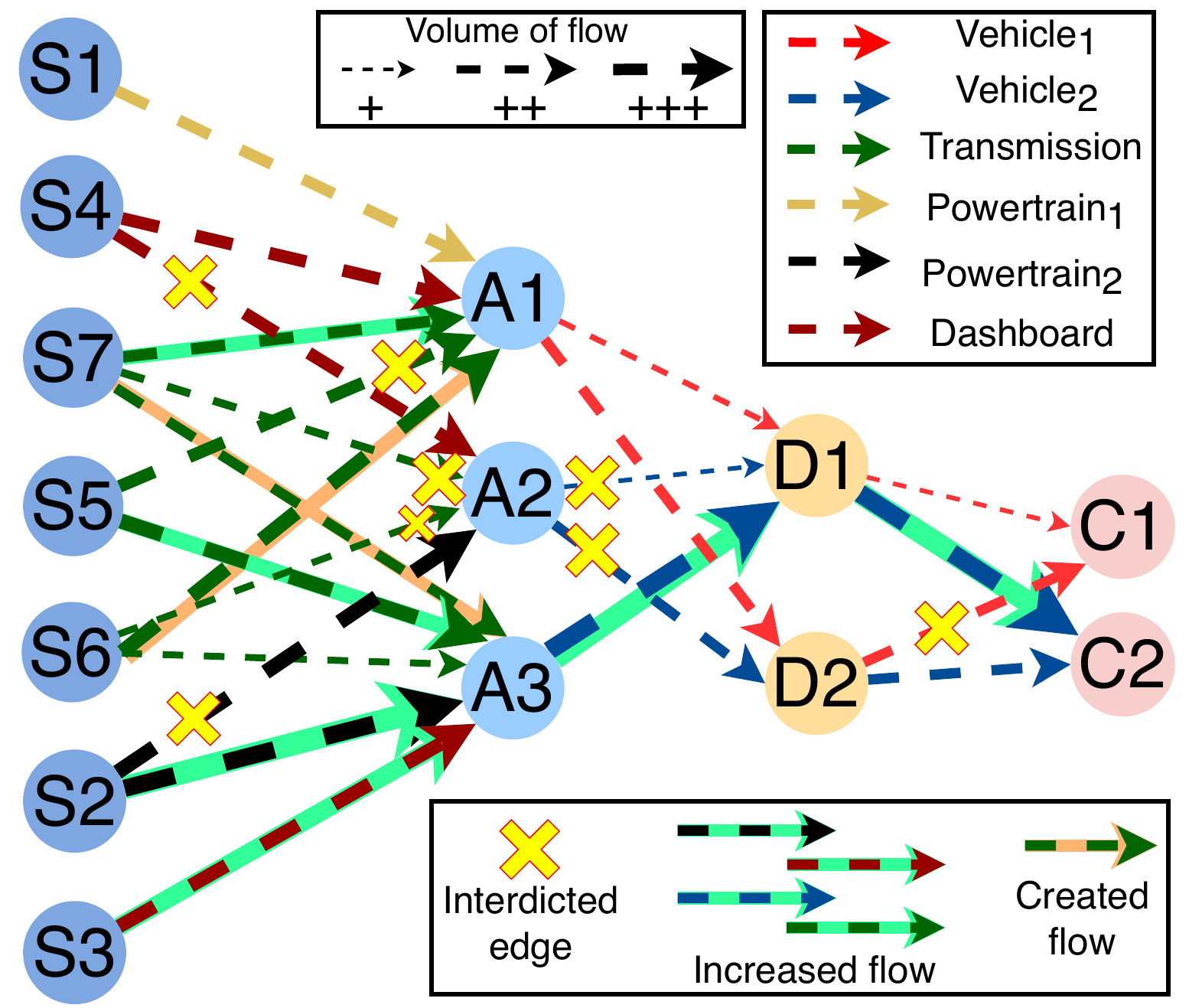}%
\label{fig:flows_A2_rat_5000}}
\caption{Optimal flows from MIP model \eqref{mip} given a disruption of vertex A2 for different lateness policies.}
\label{fig:A2_disrp_solutionflows}
\end{figure*}

\begin{figure*}[t]
\centering
\subfigure[Optimal plan for baseline tree structure.]{\includegraphics[width=0.6\columnwidth]{figures/baseline_solutionflows.pdf}%
\label{fig:Tree}}
\hfil
\subfigure[Optimal plan for reverse tree structure.]{\includegraphics[width=0.6\columnwidth]{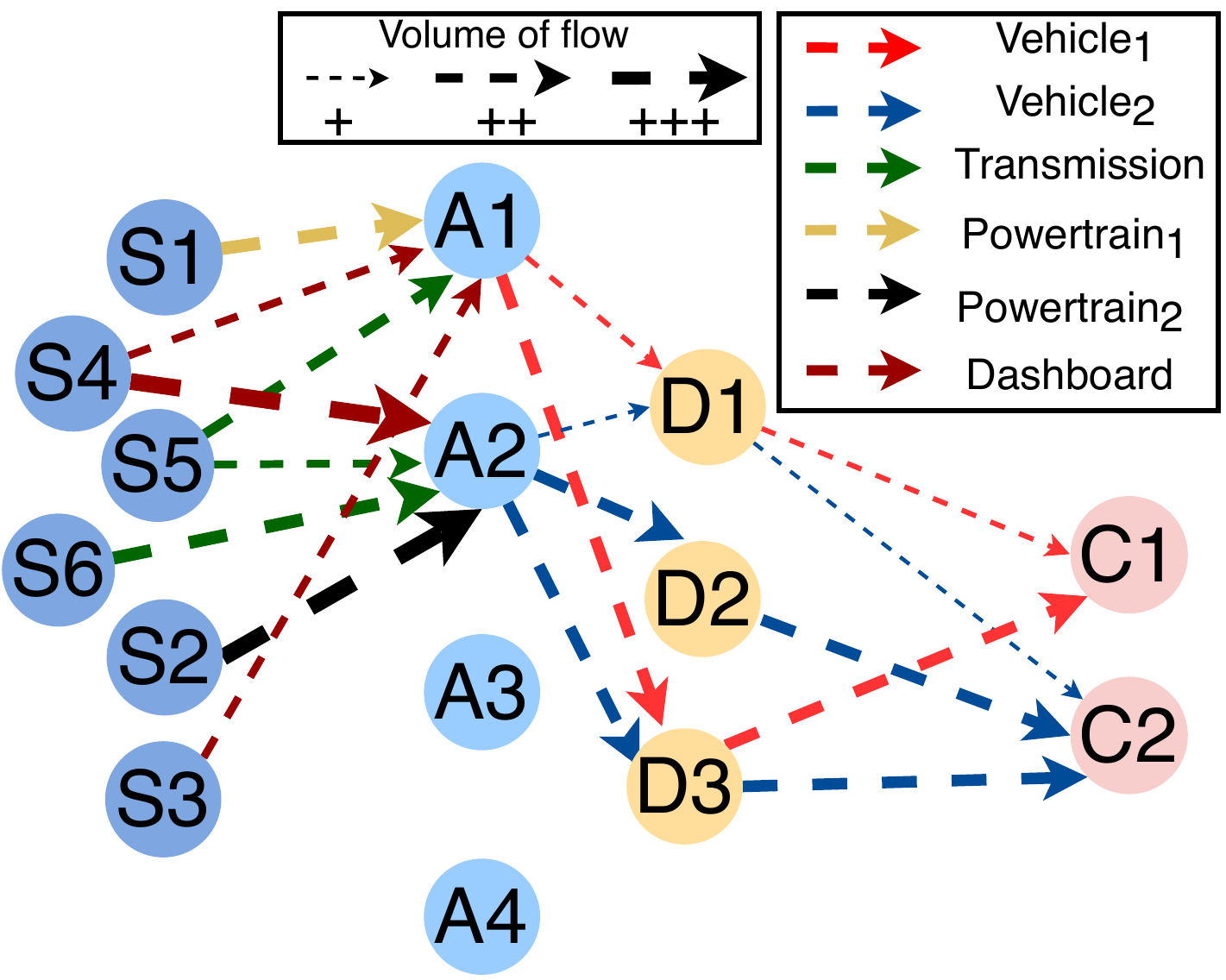}%
\label{fig:Reverse}}
\hfil
\subfigure[Optimal plan for chain structure.]{\includegraphics[width=0.6\columnwidth]{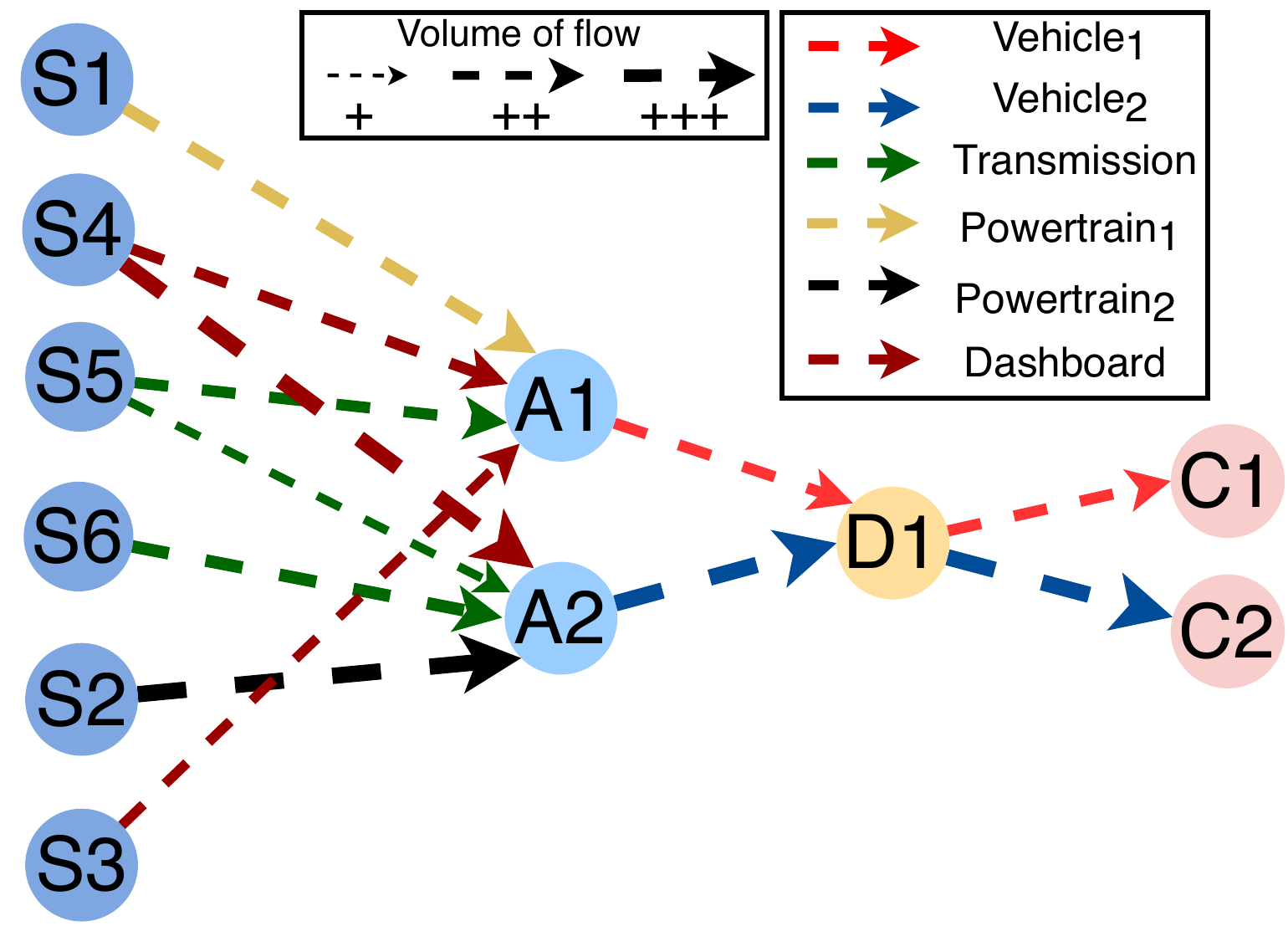}%
\label{fig:chain}}
\caption{Different SC structures to satisfy the demand of customers.}
\label{fig:othertopologies}
\end{figure*}

We perform numerical studies for the instances we describe in Section \ref{sec:outofsampletest}, on which we test the performance of solutions to Model~\eqref{mip}. For each numerical instance, we first solve Model \eqref{mip} and obtain the optimal flows $y_{ijk}$. Then, we run 300 replications of the simulation framework in which we fix $y_{ijk}$ to test average performance given stochastic lead times $l_{ijk}$. We compare simulated performance of each solution, given different parameter configuration (disruption depth, lateness policy, and SC network structure). We compute delivery lateness of final products to end customers as the performance metric. This metric can be seen as a general metric of SC effectiveness to provide timely supply.

\begin{figure}[tb]
\centerline{\includegraphics[width=0.8\columnwidth]{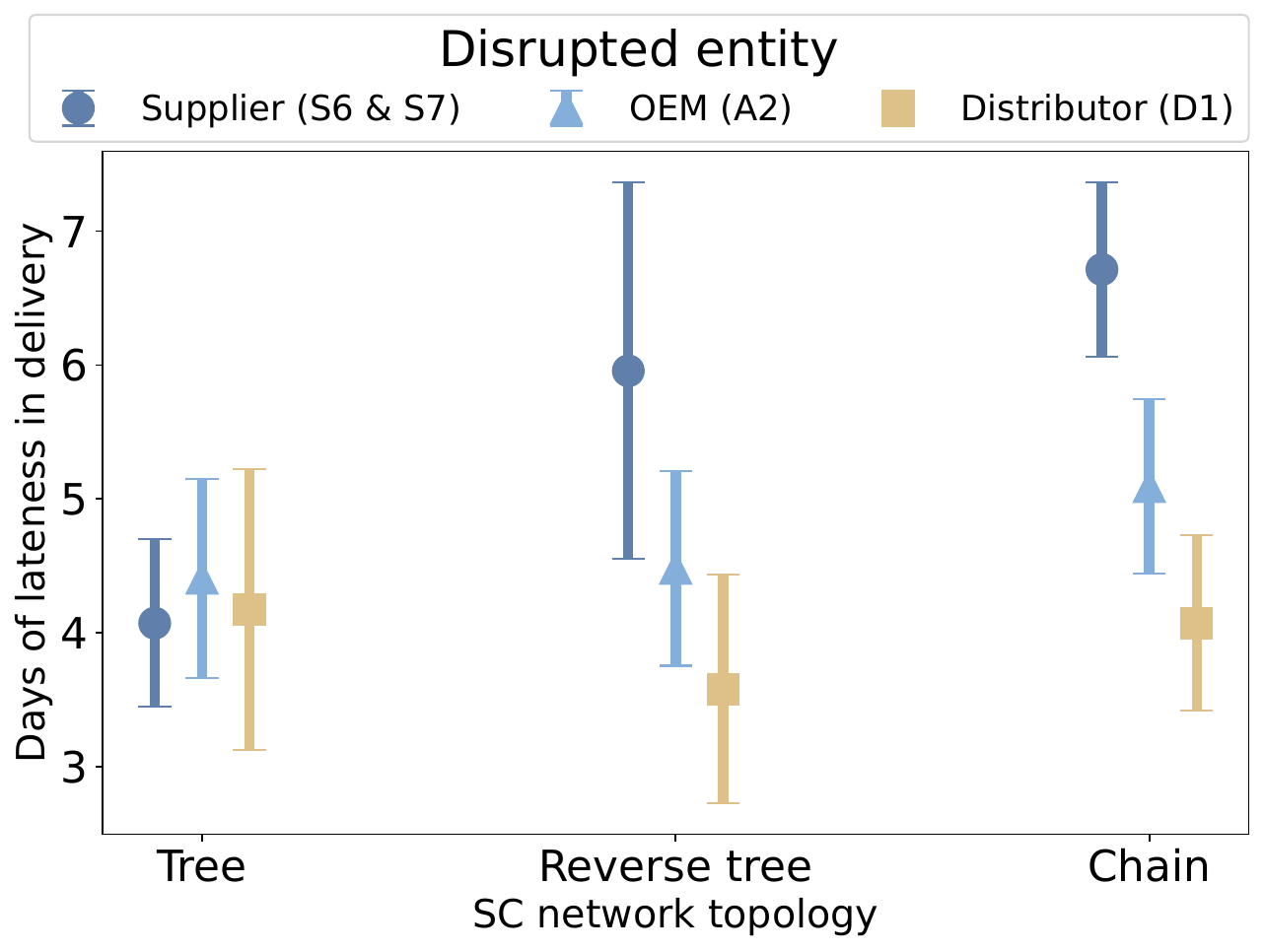}}
\caption{Average delivery lateness of different SC network typologies under diverse entity disruption.}
\vspace{-10pt}
\label{fig:topology_performance}
\end{figure}

\subsection{Case study results}
\label{sec:results}

\subsubsection{Solution of baseline}
\label{sec:baseline_solution}

We solve the MIP model~\eqref{mip} with the parameter design described in Section~\ref{sec:baselineScenario}. Fig.\ \ref{fig:baseline_flows} shows the flow solutions, which we consider as the \emph{optimal} SC operational  plan. 

\subsubsection{SC network redesign for different lateness policies}
\label{sec:latenesspolicies}

We consider the disruption scenarios described in Section~\ref{sec:SCentityDisrupts} and apply the different lateness penalties described in Section~\ref{sec:ArrLatePoli} to our proposed MIP model~\eqref{mip} to obtain optimal solutions for each scenario. 
We define three categories to describe how the network changes compared to the baseline scenario in the order of SC redesign severity. (i) K: Keep current SC, flow volumes, and incur the lateness penalty. (ii) E: Interdict edges connecting the disrupted entity with downstream entities and redesign flow solutions. (iii) V: Interdict entity by redesigning the SC without the disrupted entity. Table~\ref{tab:TopologyResults} shows that given the hierarchical structure of the SC network~\cite{kim2015supply}, as the depth of the disrupted agent increases, SC redesign is encouraged. This explains why in the case of Fig.\ \ref{fig:S7} most penalties on lateness cause the interdiction of the entity, and by the contrary, in Fig.\ \ref{fig:D1}, only with a very high penalty, the SC layout changes.
\input{tables/results}

Fig.\ \ref{fig:resultsb_simulation} demonstrates the effectiveness of the strategies prescribed by our lead-time-aware MIP model, to mitigate the effect of disruptions in the delivery of final products. From the distribution of days of lateness, MIP model~\eqref{mip} solutions enjoy less/shorter lateness, when we simulate the uncertain lead time in the SC network.

Fig.\ \ref{fig:A2_disrp_solutionflows} shows the three discussed categories for the disruption of entity A2. Since A2 has a depth of two, only when the fixed penalty is large enough, will the SC redesign reduce its influence as much as possible.

\subsubsection{Network structure in SC performance and disruption response}
\label{sec:topologystudy}

 We consider additional SC network structures described in Section~\ref{sec:SCtopologyIntro}. 
 %solving the MIP model~\eqref{mip}, such that we can test disruption response. 
 Fig.\ \ref{fig:othertopologies} shows the optimal flows for each topology without disruptions. The reverse tree structure in Fig.\ \ref{fig:Reverse} makes use of the additional distributor $D3$, but keeps as backups OEM $A3$ and the additional OEM $A4$. The chain structure has  lower cost in comparison to the other two structures. Fig.\ \eqref{fig:chain} shows the caveat of this topology, as it relies on most entities being used at full capacity. 

We perform disruption response studies similar to those in Section \ref{sec:latenesspolicies}. Fig.\ \ref{fig:topology_performance} illustrates the average lateness and standard deviation as error bars of each network structure over all of the lateness policies. We observe that the tree structure, which establishes backups further back in the SC, shows higher reliability, with a consistent performance invariant to the disrupted tier of the SC. The reverse tree structure has a better response to a disruption closer to the end customers, but presents high variability on the disruption of entities further back in the SC. This SC topology is shown to be less resilient against disruptions, even if it has multiple backups for some tiers of the SC. Finally, the chain structure performance suffers the most under entity disruptions. Only a limited number of changes to the structure can take place, making almost all entities critical.

In addition to the response strategies discussed in Section \ref{sec:latenesspolicies}, we identify an additional action shown in the our MIP model~\eqref{mip}. When considering a policy without fixed penalties (ratio 1:0), the optimal solution prescribes strategy (R) to reduce the volume of affected flows while keeping the SC layout. We interpret this as a last-resort strategy, when it is inevitable to incur in unmet demand penalties, as a trade-off to use the released capacity to work on other products that still can be delivered on time. Table~\ref{tab:TopologyResults} shows the type of responses we identify for the different scenarios. Results in Table~\ref{tab:TopologyResults} demonstrate that the chain structure has the least flexibility, as in most cases only the strategy (K) to absorb lateness penalties is possible. Finally, only the tree structure avoids strategy (R), while meeting demand in all scenarios.

\subsection{Managerial insights}
\label{sec:insight}

From the previous results shown in Section \ref{sec:casestudy}, we point out the following managerial insights. 

\begin{itemize}
    \item Lead time is a critical factor. The temporal component of material flows can make an optimal solution with the least cost be highly undesirable because of lead-time disruption and delivery lateness.

    \item A trade-off exists between disruption response and mitigation strategies, with fixed penalty policies favoring costly network-wide redesign to keep flows on time, and unit penalty policies mitigating disruptions gradually with local modifications of the network structure.

    \item Network-wide changes include redesigning a network by deleting disrupted entities. Local modifications include minimizing the influence of disrupted entities by reducing flows through them, and redistributing them to other available entities. 

    \item As the effect of lead-time disruption of entities accumulates over the sequence of flows in the SC network, entity redundancy further back into the SC (e.g., suppliers) yields more reliability. In contrast, redundancy of entities closer to end customers (e.g., distributors) shows less resiliency, since disruptions accumulate less over the remaining sequence of flows, having a lower influence on overall SC performance.
\end{itemize}

%% Conclusions 

\section{Conclusion}
\label{sec:conclusion}

In this paper, we developed an MIP model to track the delivery time of flows throughout an SC network, such that one can incorporate lateness penalties into the objective. Via numerical studies, we show that a hierarchical network structure with backups in the initial suppliers yields consistent SC performance given lead-time uncertainty of the SC entities, as opposed to having backups further down the SC.

%{\color{purple} We note that by employing a convex combination of the different elements of the objective function could disregard Pareto-optimal solutions, and for larger instances, for which we recommend tuning and a robust out-of-sample testing to evaluate the optimal flows. }

%Through the design of experiments with lead-time disruptions and using simulation, we consider the stochastic nature of lead time for solution evaluation. 
In future research, one can consider communication and coordination time-lags to incorporate real-life environments, as well as stochastic lead time and demand.  To solve a stochastic MIP formulation, we will explore the use of decomposition and cutting-plane approaches for large-scale SCs. We will also explore decentralized optimization and control methods for SC disruption response. Finally, we can explore Pareto-optimal solutions using more sophisticated multi-objective optimization approaches.

%\bibliographystyle{IEEEtran.bst}
%\bibliography{refs.bib}

\end{document}

%% file: tables/notations.tex
\begin{table}[tb]
\caption{A summary of notation for parameters and variables.}
\label{tab:notation}
    \begin{tabularx}{\columnwidth}{lX}
    \hline
    \multicolumn{2}{l}{\textbf{Input Parameters}}\\
    \hline
    $d_{ik}$& demand of product $k$ at entity $i$; by convention, $d_{ik} < 0$. \\
    $f_{ij}$& fixed transportation cost from entity $i$ to entity $j$.\\ 
    $c_{ijk}$& unit transportation cost of flowing product $k$ from entity $i$ to $j$.\\ 
    $q_{ij}$& mixed-flow capacity of edge $(i, j)$. \\ 
    $\bar{p}_{i}$& mixed-product production capacity available at entity $i$.\\
    $e_{ik}$ & production cost per unit  product $k$ at entity $i$.\\ 
    $r_{kk'}$& conversion rate from product $k$ to product $k'$, i.e., the units of product $k$ consumed to produce one unit successor product $k'$.\\ 
    $\phi_{i}$& fixed cost of opening a production line at entity $i$. \\
    $I_{ik}^0$& initial inventory of product $k$ at entity $i$ at each period.  \\
    $h_{ik}$& unit holding cost of product $k$ at entity $i$.\\ 
    $\rho_{ik}^d$& penalty per unit of unsatisfactory of demand $d_{ik}$. \\
    $t_{ik}$ & time at which entity $i$ requires the demand of product $k$. \\
    $l_{ijk}$ & lead time of product $k$ flowing from entity $i$ to entity $j$. \\
    $c^{\text{f}}_{ijk}$ &fixed late-delivery penalty of product $k$ flowing from entity $i$ to entity $j$. \\
    $c^{\text{u}}_{ijk}$ & penalty per unit of late delivery of flow  of product $k$  from entity $i$ to entity $j$. 
    \vspace{0.3em}
    \\
    \hline
    \multicolumn{2}{l}{\textbf{Decision Variables}}\\
    \hline
    $y_{ijk}$& units of product $k$ flowing from entity $i$ to entity $j$.   \\
    $\beta_{ijk}$& binary variable, equal to 1 if edge $(i, j) \in E$ is used to transport product $k$, and 0 otherwise.\\
    $x_{ik}$& actual satisfied demand of product $k$ at entity $i$.\\ 
    $p_{ik}$& units of product $k$ produced at entity $i$.\\
    $\zeta_{i}$& binary variable, equal to 1 if entity $i$ produces/assembles materials, and 0 otherwise. \\
    $I_{ik}$& inventory of product $k$ at entity $i$ by the end of time period.\\ 
    $\Delta_{ik}^d$& units of unsatisfied demand of product $k$ at entity $i$.\\
    $a_{ijk}$ & time of delivery of product $k$ to entity $j$ from entity $i$.\\
    $o_{ik}$ & time at which entity $i$ can process product $k$.\\
    $w_{ijk}$ & total penalty of late delivery of product $k$ from entity $i$ to entity $j$.\\
    $z_{ijk}$ & binary variable, equal to 1 if flow from entity $i$ to entity $j$ of product $k$ is delivered late, and 0 otherwise.
    \vspace{0.3em}
    \\
    \hline
    \end{tabularx}
\end{table}

%% file: tables/results.tex
\begin{table}[t]
\caption{Identified strategy from MIP model~\eqref{mip} optimal solution to entity disruptions under different network structures and lateness policies}
\centering
%\begin{tabularx}{\columnwidth}{cc@{\hskip 3em}c@{\hskip 2em}cccc}
\begin{tabularx}{\columnwidth}{c@{\hskip 3em}c@{\hskip 1em}c@{\hskip 1em}c@{\hskip 1em}c@{\hskip 1em}c@{\hskip 1em}c@{\hskip 1em}c@{\hskip 1em}c@{\hskip 1em}c@{\hskip 1em}}
\hline
\multicolumn{3}{c}{\textbf{Scenario}}& \multicolumn{4}{c}{\textbf{(Unit penalty : Fixed penalty)}}\vspace{0.2em}\\
Agent & Depth & Structure & No penalty & 1:0 & 1:500 & 1:5000 \vspace{0.3em}\\
\hline 
\noalign{\vskip 0.3em}
\multirow{3}{1em}{Supplier (S7, S6)}& \multirow{3}{1em}{3} &Tree &K &V& V & V \vspace{0.3em}\\
& &Reverse Tree &K & R & K & K \vspace{0.3em}\\
& &Chain &K & R & K & K \vspace{0.3em}\\
\hline
\multirow{3}{1em}{OEM (A2)}& \multirow{3}{1em}{2} & Tree &K &K& E & V \vspace{0.3em}\\
& & Reverse Tree &K &K& E & V \vspace{0.3em}\\
& & Chain &K &K& K & K \vspace{0.3em}\\
\hline
\multirow{3}{1em}{Distributor (D1)}&\multirow{3}{1em}{1}& Tree &K &K& K & E \vspace{0.3em}\\
& & Reverse Tree &K &V& V & V \vspace{0.3em}\\
& & Chain &K &K& K & K \vspace{0.3em}\\
\hline
\multicolumn{7}{l}{\footnotesize K: keep flows; E: interdict edges; V: interdict entity; R: reduce flows.} \\
\end{tabularx}
\label{tab:TopologyResults}
\end{table}